\documentclass[letterpaper, 10 pt, conference]{ieeeconf}  

\IEEEoverridecommandlockouts                              

\usepackage{graphics} 
\usepackage{epsfig} 
\usepackage{mathptmx} 
\usepackage{amsmath} 
\usepackage{amssymb}  

\usepackage{mathrsfs}  

\usepackage[dvipsnames]{xcolor}

\usepackage{tikz,pgfplots}

\usepackage{comment}

\newtheorem{defn}{Definition}[section]
\newtheorem{thm}{Theorem}[section]
\newtheorem{remark}{Remark}[section]

\newtheorem{lemma}{Lemma}[section]
\newtheorem{cor}{Corollary}[section]

\renewcommand{\d}{\mathrm{d}}
\renewcommand{\L}{\mathrm{L}}
\renewcommand{\H}{\mathrm{H}}

\title{\LARGE \bf
Spectral analysis of a class of linear hyperbolic partial differential equations
}

\author{Anthony Hastir$^{1}$, Birgit Jacob$^{1}$ and Hans Zwart$^{2}$
\thanks{*This research was conducted with the financial support of the German Research Foundation (DFG). A. H. is supported by the DFG via the Grant HA 10262/2-1.}
\thanks{$^{1}$Anthony Hastir and Birgit Jacob are with the University of Wuppertal, School of Mathematics and Natural Sciences, Gaußstr. 20, 42119 Wuppertal, Germany,
        {\tt\small hastir@uni-wuppertal.de, bjacob@uni-wuppertal.de}}%
\thanks{$^{2}$Hans Zwart is with the Department of Applied Mathematics, University of Twente, 7500 AE Enschede, The Netherlands, and Department of Mechanical Engineering, Technische Universiteit Eindhoven, 5600 MB Eindhoven,
The Netherlands,
        {\tt\small h.j.zwart@utwente.nl}}%
}

\begin{document}

\maketitle
\thispagestyle{empty}
\pagestyle{empty}

\begin{abstract}
A class of linear hyperbolic partial differential equations, sometimes called networks of waves, is considered. For this class of systems, necessary and sufficient conditions are formulated on the system matrices for the operator dynamics to be a Riesz-spectral operator. In that case, its spectrum is computed explicitly, together with the corresponding eigenfunctions, which constitutes the main result of our note. In particular, this enables to characterize easily many different concepts, such as stability. We apply our results to characterize exponential stability of a co-current heat exchanger.
\end{abstract}

\begin{keywords}
Discrete Riesz-spectral operators, Distributed parameter systems, Hyperbolic partial differential equations, Stability of linear systems
\end{keywords}


\section{INTRODUCTION}
Hyperbolic partial differential equations (PDEs) constitute a large class of dynamical systems. The references \cite{Hyperbolic_Phillips_1957} and \cite{Hyperbolic_Phillips_1959} already dealt with linear hyperbolic PDEs, studied in the framework of semigroup theory. From an application point of view, hyperbolic PDEs are able to describe many different physical phenomena. In their nonlinear formulation, they can model the evolution of flow of water and sediments in open channels, via the Saint-Venant-Exner equations. As another example, the dynamics of plug-flow reactors are also governed by nonlinear hyperbolic PDEs. Furthermore, hyperbolic PDEs are able to describe the motion of an inviscid ideal gas in a rigid cylindrical pipe, via the nonlinear Euler equations. Once linearized, they can model the propagation of current and voltage in transmission lines, or the dynamics of vibrating strings to mention a few. More generally, conservation laws may be written via hyperbolic PDEs, both linear and nonlinear. Beyond these examples, a comprehensive overview of many applications governed by hyperbolic PDEs may be found in \cite{BastinCoron_Book}. Therein, these equations are presented both in their nonlinear and linear forms, where aspects like stability and control are described as well. 

An important property that some linear infinite-dimensional systems may exhibit is the Riesz-spectral feature of their operator dynamics. In addition, if the resolvent operator associated to the operator dynamics is compact, the spectrum of such operators is composed only of a discrete set of eigenvalues, which turns out to be useful to study properties like stability for instance. Historically the concept of spectral operator has been introduced in \cite{Dunford_PartIII}. Few years later, discrete spectral operators were discussed in \cite{Curtain_Riesz} while Riesz-spectral operators came in \cite{CurtainZwart1995}, in which the eigenvalues are assumed to be simple. A generalization to the multiple case may be found in \cite{GuoZwart_2001Riesz}. Many different works investigating the Riesz-spectral property of linear operators describing hyperbolic PDEs may be found in the literature. For instance, \cite{Guo_Riesz_2003} studied this question for coupled string equations whereas \cite{Guo_SIAM_2001} and \cite{Guo_SIAM_2002} are dealing with Riesz-spectral property for beam equations. We refer the reader to the book \cite{Guo2019_Book} for an overview of different applications in which Riesz-spectral operators may be encountered. Generators of $C_0$-groups and Riesz-spectral operators have been shown to be related in some sense, as the works \cite{Guo_2004_SCL} and \cite{Zwart_2010} demonstrate.  A few years later, port-Hamiltonian systems of the form described in \cite{JacobZwart} have been shown Riesz-spectral systems if and only if the operator dynamics generates a $C_0$-group, see \cite{JacobKaiserZwart}.

In this work, we consider the following linear hyperbolic PDEs
\begin{align}
\frac{\partial \tilde{z}}{\partial t}(\zeta,t) &= -\frac{\partial}{\partial \zeta}(\lambda_0(\zeta)\tilde{z}(\zeta,t)) + M(\zeta)\tilde{z}(\zeta,t),\nonumber\\
\tilde{z}(\zeta,0) &= \tilde{z}_0(\zeta),\label{StateEquation_Diag_Uniform}
\end{align}
where $t\geq 0$ and $\zeta\in [0,1]$ are the time and the space variables, respectively. The function $\tilde{z}(\zeta,t)\in\mathbb{C}^n$ is the state with initial condition $\tilde{z}_0\in\L^2(0,1;\mathbb{C}^n) =: X$. The PDEs are complemented with the following boundary conditions
\begin{align}
0 &= \lambda_0(0) K\tilde{z}(0,t) + \lambda_0(1) L \tilde{z}(1,t),\label{BC_Diag_Uniform}
\end{align}
$t\geq 0$. The scalar-valued function $\lambda_0\in \L^\infty(0,1;\mathbb{R}^+)$ satisfies $\lambda_0(\zeta)\geq \epsilon > 0$ a.e.~on $[0,1]$. $M(\zeta)$ is a space-dependent square matrix that satisfies $M\in\L^\infty(0,1;\mathbb{C}^{n\times n})$ while $K\in\mathbb{C}^{n\times n}$ and $L\in\mathbb{C}^{n\times n}$. 

This class of hyperbolic systems appears when looking at networks of conservation laws such as networks of electric lines, chains of density-velocity systems or genetic regulatory networks to mention a few, see e.g. \cite{BastinCoron_Book} for an overview of the covered applications. We also refer the reader to the book \cite{Luo_Guo} in which a more general version of the class \eqref{StateEquation_Diag_Uniform}--\eqref{BC_Diag_Uniform} is presented and analyzed. 

The main result of this note consists in characterizing the Riesz-spectral property of the operator dynamics associated to \eqref{StateEquation_Diag_Uniform}--\eqref{BC_Diag_Uniform}, where Riesz-spectral has to be understood in the sense of \cite[Definition 2.2]{JacobKaiserZwart}. Moreover, we give explicit analytical expressions of the eigenvalues and the eigenfunctions of that operator, that can be helful in showing several dynamical properties. This characterization is used to describe exponential stability of a co-current heat exchanger with specific boundary conditions.

\section{WELL-POSEDNESS}\label{Section_Well-Posed}

The PDE \eqref{StateEquation_Diag_Uniform} with the boundary conditions \eqref{BC_Diag_Uniform} may be written in an abstract way as $\dot{\tilde{z}}(t) = A\tilde{z}(t), \tilde{z}(0) = \tilde{z}_0$ with the (unbounded) linear operator $A$ defined as
\begin{align*}
Af &= -\frac{\d}{\d\zeta}(\lambda_0 f) + Mf,\\
D(A) &= \{f\in X \mid \lambda_0 f\in\H^1(0,1;\mathbb{C}^n),\\
&\hspace{2.5cm} \lambda_0(0)Kf(0) + \lambda_0(1)Lf(1) = 0\}.
\end{align*}
Here the state space $X$ is equipped with the inner product $\langle f,g\rangle_X = \displaystyle\int_0^1 \lambda_0(\zeta)g(\zeta)^*f(\zeta)\d\zeta$. 

By well-posedness of \eqref{StateEquation_Diag_Uniform}--\eqref{BC_Diag_Uniform}, we mean that the operator $A$ is the infinitesimal generator of a $C_0$-semigroup of bounded linear operators. To look at this question, observe first that the assumption $M(\cdot)\in\L^\infty(0,1;\mathbb{C}^{n\times n})$ implies that the matrix differential equation 
\begin{align}
P'(\zeta) &= \lambda_0(\zeta)^{-1}M(\zeta)P(\zeta),\hspace{1cm} P(0) = I,
\label{ODE_Matrix}
\end{align}
has a unique solution that is invertible\footnote{According to \cite[Lemma 2.1]{HastirJacobZwart_LQ}, the inverse of $P$ satisfies the following differential equation 
\begin{align*}
(P^{-1})'(\zeta) &= -\lambda_0(\zeta)^{-1}P^{-1}(\zeta)M(\zeta),\hspace{1cm} P^{-1}(0) = I.
\end{align*}
} on $[0,1]$, see \cite{trostorff2023characterisation}. By defining the multiplication operator $\mathfrak{Q}: X\to X$ by $(\mathfrak{Q}\tilde{z})(\zeta) = P^{-1}(\zeta)\tilde{z}(\zeta)$, there holds
\begin{align}
\mathfrak{Q}A\mathfrak{Q}^{-1}g &= -\frac{\d}{\d\zeta}(\lambda_0 g),\label{Op_A}\\
D(\mathfrak{Q}A\mathfrak{Q}^{-1}) &= \{g\in X \mid \lambda_0 g\in\H^1(0,1;\mathbb{C}^n),\nonumber\\
&\hspace{0.5cm} \lambda_0(0)Kg(0) + \lambda_0(1)LP(1)g(1) = 0\},\label{D(A)}
\end{align}
see \cite[Lemma 2.2]{HastirJacobZwart_LQ}. This implies that the change of variables $z(\zeta,t) = P^{-1}(\zeta)\tilde{z}(\zeta,t)$ transforms \eqref{StateEquation_Diag_Uniform}--\eqref{BC_Diag_Uniform} equivalently to the following PDE
\begin{align}
\frac{\partial z}{\partial t}(\zeta,t) &= -\frac{\partial}{\partial \zeta}(\lambda_0(\zeta)z(\zeta,t)),\label{StateEquation_Diag_Uniform_Chg_Var}\\
z(\zeta,0) &= z_0(\zeta),\hspace{0.5cm} \zeta\in (0,1),\label{InitEquation_Diag_Uniform_Chg_Var}\\
0 &= \lambda_0(0) Kz(0,t) + \lambda_0(1) L P(1)z(1,t).\label{BC_Diag_Uniform_Chg_Var}
\end{align}
Let us denote by $\mathfrak{A}$ the operator $\mathfrak{Q}A\mathfrak{Q}^{-1}$. It is easy to see that the operator $A$ is the infinitesimal generator of a $C_0$-semigroup if and only if so is $\mathfrak{A}$. Well-posedness of \eqref{StateEquation_Diag_Uniform_Chg_Var}--\eqref{BC_Diag_Uniform_Chg_Var} is characterized in the next lemma, in which the port-Hamiltonian formalism of \cite{JacobZwart} is used, see also \cite[Proposition 2.1]{HastirJacobZwart_LQ}.
\begin{lemma}\label{Lemma_Well-Posed}
The operator $\mathfrak{A}$ defined in \eqref{Op_A}--\eqref{D(A)} is the infinitesimal generator of a $C_0$-semigroup if and only if $K$ is invertible.
\end{lemma}

\proof
We rely on \cite[Theorem 1.5]{JacobMorrisZwart}. Observe that \eqref{StateEquation_Diag_Uniform_Chg_Var}--\eqref{BC_Diag_Uniform_Chg_Var} may be written as the following port-Hamiltonian system
\begin{equation}
\left\{\begin{array}{l}
\displaystyle\frac{\partial z}{\partial t}(\zeta,t) = P_1\frac{\partial}{\partial\zeta}(\mathcal{H}(\zeta)z(\zeta,t)),\hspace{0.8cm} \zeta\in (0,1),t\geq 0,\\
z(\zeta,0) = z_0(\zeta),\hspace{0.8cm} \zeta\in (0,1)\\
0 = \tilde{W}_{B}\left(\begin{matrix}
(\mathcal{H}z)(1,t)\\
(\mathcal{H}z)(0,t)
\end{matrix}\right), \hspace{0.8cm} t\geq 0,
\end{array}\right.
\label{Class_PHS}
\end{equation}
where $P_1 = -I, \mathcal{H}(\zeta) = \lambda_0(\zeta)$ and $\tilde{W}_{B} := \left(\begin{matrix}W_1 & W_0\end{matrix}\right)$ with $W_1 = LP(1)$ and $W_0 = K$. Following the notations of \cite{JacobMorrisZwart}, we denote by $Z^+(\zeta)$ the span of eigenvectors of $P_1\mathcal{H}(\zeta)$ corresponding to positive eigenvalues and $Z^-(\zeta)$ the span of eigenvectors corresponding to negative eigenvalues. We have $Z^+(\zeta) = \{0\}$ and $Z^-(\zeta) = \mathbb{C}^n$ because $\lambda_0$ is a positive function. According to \cite[Theorem 1.5]{JacobMorrisZwart} the operator $\mathfrak{A}$ is the infinitesimal generator of a $C_0$-semigroup if and only if $W_1\mathcal{H}(1)Z^+(1)\oplus W_0\mathcal{H}(0)Z^-(0) = \mathbb{C}^n$. Equivalently, if and only if $\mathbb{C}^n = L\lambda_0(1)P(1)\{0\}\oplus K\lambda_0(0)\mathbb{C}^n = K\lambda_0(0)\mathbb{C}^n$. Since $\lambda_0(0)$ is not zero, this holds if and only if $K$ is invertible.
\endproof

In order to be able to describe the spectrum of the operator $\mathfrak{A}$ only with eigenvalues, we need to see whether the semigroup generated by $\mathfrak{A}$ is also a group and whether the resolvent operator associated to $\mathfrak{A}$ is compact, see \cite{JacobKaiserZwart}. This first property is characterized in the following lemma.

\begin{lemma}\label{Lemma_Group}
The $C_0$-semigroup generated by $\mathfrak{A}$ is a $C_0$-group if and only if $L$ is invertible.
\end{lemma}

\proof
According to \cite[Theorem 2.8]{JacobKaiserZwart}, the $C_0$-semigroup generated by $\mathfrak{A}$ is a group if and only if $W_1\mathcal{H}(1)Z^-(1)\oplus W_0\mathcal{H}(0)Z^+(0) = \mathbb{C}^n$, where the notation is the same as in the proof of Lemma \ref{Lemma_Well-Posed}. Equivalently, the previous equality is $\lambda_0(1)LP(1)\mathbb{C}^{n} = \mathbb{C}^n$. Since $\lambda_0(1)$ is positive and $P(1)$ is invertible, the semigroup generated by $\mathfrak{A}$ is a group if and only if $L$ is invertible.
\endproof

The fact that the resolvent operator associated to $\mathfrak{A}$ is compact comes from \cite[Theorem 2.26]{Villegas}.

Note that characteristic based methods could be used to solve \eqref{StateEquation_Diag_Uniform_Chg_Var}--\eqref{BC_Diag_Uniform_Chg_Var}, see e.g. \cite{Mas52}.

\section{SPECTRAL ANALYSIS}\label{Section_Spectral-Analysis}

As is highlighted in \cite[Theorem 2.8]{JacobKaiserZwart}, the operator $\mathfrak{A}$ generating a $C_0$-group is equivalent to the operator $\mathfrak{A}$ being a discrete Riesz-spectral operator. Before recalling the definition of a Riesz-spectral operator, we need to introduce the concept of spectral projection, see \cite[Definition 2.1]{JacobKaiserZwart}.
\begin{defn}
For an operator $\mathfrak{A}$ on $X$ we call $\gamma\subset\sigma(A)$ a compact spectral set if $\gamma$ is a compact subset of $\mathbb{C}$ which is open and closed in $\sigma(A)$. The spectral projection on the spectral subset $\gamma$ is defined as
\begin{equation*}
E(\gamma) = \frac{1}{2\pi i}\int_\Gamma (sI - \mathfrak{A})^{-1}\d s,
\end{equation*}
where $\Gamma$ is a closed Jordan curve containing every point of $\gamma$ and no point of $\sigma(\mathfrak{A})\setminus\gamma$.
\end{defn}

Discrete Riesz-spectral operators are introduced in the next definition, see \cite[Definition 2.2]{JacobKaiserZwart}.

\begin{defn}\label{Def_Riesz}
An operator $\mathfrak{A}$ with compact resolvent and $\sigma(\mathfrak{A}) = \left(s_n\right)_{n\in\mathbb{I}}, \mathbb{I}$ countable, is a discrete Riesz-spectral operator if the following hold:
\begin{itemize}
\item For every $n\in\mathbb{I}$ there exists $N_n\in\mathcal{L}(X)$ such that $AE_n = (s_n + N_n)E_n$.
\item The sequence of closed subspaces $\left(E_n(X)\right)_{n\in\mathbb{I}}$ is a Riesz basis of subspaces of $X$, that is, $\text{span}\left(E_n(X)\right)_{n\in\mathbb{I}}$ is dense and there exists an isomorphism $T\in\mathcal{L}(X)$ such that $(T E_n(X))_{n\in\mathbb{I}}$ is a system of pairwise orthogonal subspaces of $X$.
\item $N := \sum_{n\in\mathbb{I}}N_n$ is bounded and quasi nilpotent.
\end{itemize}
\end{defn}

\begin{remark}
\begin{itemize}
\item As is mentioned in \cite{JacobKaiserZwart}, for an operator $\mathfrak{A}$ with compact resolvent and spectrum denoted by $\left(s_n\right)_{n\in\mathbb{I}}$, the notation $E_n$ stands for $E((s_n)), n\in\mathbb{I}$.
\item The eigenvalues of a Riesz-spectral operator described in Definition \ref{Def_Riesz} are not supposed to be simple, i.e. they could have algebraic multiplicity larger than $1$. Moreover, equality between the algebraic multiplicity and the geometric multiplicity of an eigenvalue is not assumed in Definition \ref{Def_Riesz}, which could lead to Jordan blocks and generalized eigenvectors. This makes Riesz-spectral operators in Definition \ref{Def_Riesz} different from those in \cite[Definition 3.2.6]{CurtainZwart2020}.
\end{itemize}
\end{remark}

We are now in place to state our main theorem, that gives a characterization for the operator $\mathfrak{A}$ to be a discrete Riesz-spectral operator. Moreover, analytical expressions of the eigenvalues and the eigenfunctions of the operator $\mathfrak{A}$ are given.

\begin{thm}\label{Main_Thm}
Consider the operator $\mathfrak{A}$ defined in \eqref{Op_A}--\eqref{D(A)}. Let $A_d := -K^{-1}LP(1)$ and $\eta(\zeta) := \int_0^\zeta\lambda_0(\tau)^{-1}\d\tau$. Let $\rho_k := \vert\rho_k\vert e^{i\theta_k},~ k = 1,\dots,n, ~\theta_k\in [0,2\pi[$, be the eigenvalues of $A_d$. The operator $\mathfrak{A}$ is a discrete Riesz-spectral operator in the sense of Definition \ref{Def_Riesz} if and only if the matrices $K$ and $L$ are invertible. In this case, the following hold
\begin{itemize}
\item The eigenvalues of $\mathfrak{A}$ are given by
\begin{equation}
\mu_{kl} := \eta(1)^{-1}\left[\log(\vert\rho_k\vert) + i\left(\theta_k + 2\pi l\right)\right],
\label{Eigenvalues_A}
\end{equation}
$k = 1,\dots, n$, and $l\in\mathbb{Z}$;
\item If $0\neq v_k\in\mathrm{ker}(\rho_k-A_d)$ is an eigenvector of $A_d$, then the eigenvector of $\mathfrak{A}$ corresponding to $\mu_{kl}$ is given by 
\begin{align}
\phi_{kl}(\zeta) = \lambda_0(0)\lambda_0(\zeta)^{-1}e^{-\mu_{kl}\eta(\zeta)}v_k.
\label{EigenFct_Order1}
\end{align}
\item Let $\{(\rho_k-A_d)^jv\}_{j=0}^{p-1}$ for some $v\in\mathbb{C}^n$ and some $p\in\mathbb{N}$ be a cycle of generalized eigenvectors of $A_d$. Moreover, let the sequence of non-zero vectors $\{\omega_m\}_{m=1}^{p}$ be constructed such that $\omega_m$ is in the span of $\{(\rho_k-A_d)^{p-r}v\}_{r=1}^{m}$ and
\begin{align}
(\rho_k-A_d)\omega_1 &= 0,\nonumber\\
(\rho_k-A_d)\omega_m &= A_d\sum_{q=1}^{m-1}\Omega_q(1)\omega_{m-q},\label{Construction_omega}
\end{align}
where $\Omega_m$ is defined by induction as 
\begin{align*}
    \Omega_m(\zeta) &= \int_0^{\zeta}\lambda_0(\alpha)^{-1}\Omega_{m-1}(\alpha)\d\alpha, \,\,\,m\geq 1,\,\,\,\Omega_0(\zeta) = 1.
\end{align*}
The non-zero vectors $\phi_{klj}\in D(\mathfrak{A}), j=1,\dots, p$ given by
\begin{align} 
&\phi_{klj}(\zeta) = \lambda_0(0)\lambda_0(\zeta)^{-1}e^{-\mu_{kl}\eta(\zeta)}\left[\omega_j+\right.\nonumber\\
&\left. \sum_{m=0}^{j-2}\omega_{j-1-m}\Omega_{m+1}(\zeta)\right]\label{Gen_EigFct}
\end{align}
are generalized eigenvectors of $\mathfrak{A}$ corresponding to $\mu_{kl}$. Moreover, $(\mu_{kl}-\mathfrak{A})^j\phi_{klj} = 0$ and $(\mu_{kl}-\mathfrak{A})^{j-1}\phi_{klj} \neq 0$.
\item In the case where $A_d$ has a basis of eigenvectors, denoted by $\{v_k\}_{1\leq k\leq n}$, the set $\{\phi_{kl}\}_{1\leq k\leq n, l\in\mathbb{Z}}$ is complete and it can be made orthogonal with respect to the following inner product
\begin{equation}
\langle f,g\rangle_\mathcal{W} := \int_0^1 g^*(\zeta)\lambda_0(\zeta)\mathcal{W}(\zeta)f(\zeta)\d\zeta,
\label{Weighted_InnerProduct}
\end{equation}
where the positive definite matrix $\mathcal{W}(\zeta)$ is defined with $v_k^*\mathcal{W}(\zeta)v_l = w_k(\zeta)\delta_{kl}$, with
\begin{equation*}
w_k(\zeta) := \displaystyle e^{2\eta(1)^{-1}\eta(\zeta)\log\vert\rho_k\vert}\vert \lambda_0(0)\vert^{-2}\eta(1)^{-1}.
\end{equation*}
\end{itemize}
\end{thm}
\vspace{0.2cm}
\proof
It is shown in \cite[Theorem 2.8]{JacobKaiserZwart} that $\mathfrak{A}$ being a discrete Riesz-spectral operator is equivalent for $\mathfrak{A}$ to be the generator of a $C_0$-group. According to Lemmas \ref{Lemma_Well-Posed} and \ref{Lemma_Group}, $\mathfrak{A}$ is a discrete Riesz-spectral operator if and only if both $K$ and $L$ are invertible. In order to establish \eqref{Eigenvalues_A} and \eqref{EigenFct_Order1}, let us focus on the following eigenvalue problem, where we seek for functions $0\neq \phi$ in $D(\mathfrak{A})$ and scalars $\mu\in\mathbb{C}$ such that
$\mathfrak{A}\phi = \mu\phi$. Equivalently, this equation may be written as
\begin{align*}
\frac{\d}{\d\zeta}(\lambda_0(\zeta)\phi(\zeta)) = -\frac{\mu}{\lambda_0(\zeta)}(\lambda_0(\zeta)\phi(\zeta)).
\end{align*}
The solution of the above differential equation is given by
\begin{equation}
\phi(\zeta) = \lambda_0(\zeta)^{-1}e^{-\mu \eta(\zeta)}\lambda_0(0)c.\label{Intermediate_phi}
\end{equation}
Now we determine $\mu$ and $c$ via the boundary conditions induced by $D(\mathfrak{A})$. The condition $\phi\in D(\mathfrak{A})$ is equivalent to 
\begin{equation}
\lambda_0(0)K\phi(0) + \lambda_0(1)LP(1)\phi(1) = 0.
\label{Intermediate_BC}
\end{equation}
Plugging the expression \eqref{Intermediate_phi} into \eqref{Intermediate_BC} yields
\begin{equation}
e^{\mu \eta(1)}c + K^{-1}LP(1)c = 0.
\label{Eigenvalue_Equation}
\end{equation}
Before continuing, note that if this equation holds for some $\mu$ and some $c$, then it also holds when $\mu$ is replaced by $\tilde{\mu} := \mu + 2 i\pi n \eta(1)^{-1}$, for all $n\in\mathbb{Z}$, with the same $c$. Let us denote by $\{\rho_k\}_{k=1,\dots,n}$ the eigenvalues of the matrix $A_d$. Because the matrices $K, L$ and $P(1)$ are invertible, so is $A_d$, which implies that none of the eigenvalues $\rho_k$ can be $0$. From \eqref{Eigenvalue_Equation}, it is clear that $c\in\mathrm{ker}(e^{\mu \eta(1)}-A_d)$ and that $e^{\mu \eta(1)}$ is an eigenvalue of $A_d$. Now let us take some $1\leq k\leq n$. For that $k$, we may define $\mu$ as
$e^{\mu \eta(1)} = \rho_k$ such that \eqref{Eigenvalue_Equation} holds. Remark that $\rho_k$ may be written as $\rho_k := \vert\rho_k\vert e^{i\theta_k}$ with $\theta_k\in [0,2\pi[$. As a consequence, all the $\mu$'s that satisfy \eqref{Eigenvalue_Equation} are given by $\mu_{kl} := \eta(1)^{-1}\left[\log(\vert\rho_k\vert) + i\left(\theta_k + 2\pi l\right)\right], k = 1, \dots, n$, and $l\in\mathbb{Z}$. We will show by induction that the construction of the non-zero vectors $\{\omega_m\}_{m=1}^p$ satisfying \eqref{Construction_omega} and such that 
\begin{equation}
\omega_m\in\textrm{span}\{(\rho_k-A_d)^{p-1}v, \dots, (\rho_k-A_d)^{p-m}v\},
\label{Property_omega_m}
\end{equation}
is possible. We start with $m=1$. We choose $\omega_1 = (\rho_k-A_d)^{p-1}v$. It is then clear that \eqref{Construction_omega} holds because of the cycle property $(\rho_k-A_d)^pv = 0$. Now assume that the construction was successful until $m=J,\,\,\, J<p$. Observe that 
\begin{align*}
&(\rho_k-A_d)\omega_{J+1}\\
&= A_d\sum_{q=1}^{J}\Omega_q(1)\omega_{J+1-q}\\
&= (A_d-\rho_k)\sum_{q=1}^{J}\Omega_q(1)\omega_{J+1-q} + \rho_k\sum_{q=1}^{J}\Omega_q(1)\omega_{J+1-q}.
\end{align*}
By the induction assumption, the first term in the above sum is in 
\begin{align*}
\textrm{span}\{(\rho_k-A_d)^{p-1}v, \dots, (\rho_k-A_d)^{p-J-1}v\}
\end{align*}
while the second term is in 
\begin{align*}
\textrm{span}\{(\rho_k-A_d)^{p-1}v, \dots, (\rho_k-A_d)^{p-J}v\}.
\end{align*}
It follows that $(\rho_k-A_d)\omega_{J+1} = \sum_{q=1}^J\beta_q(\rho_k-A_d)^{p-q}v$ for some $\{\beta_q\}_{q=1}^J$. We see that by choosing
\begin{align*}
\omega_{J+1} = \sum_{q=1}^J\beta_q(\rho_k-A_d)^{p-q-1}v,
\end{align*}
we satisfy both criteria. We will now establish \eqref{Gen_EigFct}. We start by applying the operator $\mu_{kl}-\mathfrak{A}$ to \eqref{Gen_EigFct}. There holds
\begin{align}
&\mu_{kl}\phi_{klj} - \mathfrak{A}\phi_{klj}\nonumber\\
&= \mu_{kl}\phi_{klj}(\zeta) + \frac{\d}{\d\zeta}(\lambda_0(\zeta)\phi_{klj}(\zeta))\nonumber\\
&=\lambda_0(0)e^{-\mu_{kl}\eta(\zeta)}\sum_{m=0}^{j-2}\omega_{j-1-m}\lambda_0^{-1}(\zeta)\Omega_m(\zeta)\nonumber\\
&= \lambda_0(0)\lambda_0(\zeta)^{-1}e^{-\mu_{kl}\eta(\zeta)}\left[\omega_{j-1}+\sum_{m=1}^{j-2}\omega_{j-1-m}\Omega_m(\zeta)\right]\nonumber\\
&= \lambda_0(0)\lambda_0(\zeta)^{-1}e^{-\mu_{kl}\eta(\zeta)}\left[\omega_{j-1}\right.\nonumber\\
&\left.+\sum_{\tilde{m}=0}^{(j-1)-2}\omega_{(j-1)-1-\tilde{m}}\Omega_{\tilde{m}+1}(\zeta)\right]\label{Intermediate_GenFct}
\end{align}
where $\frac{\d}{\d\zeta}\Omega_{m+1}(\zeta) = \lambda_0(\zeta)^{-1}\Omega_m(\zeta)$ has been used. It is clear that \eqref{Intermediate_GenFct} is of the same form as \eqref{Gen_EigFct} with $j$ replaced by $j-1$. Applying again the operator $\mu_{kl}-\mathfrak{A}$ will remove successively the terms in the sum in \eqref{Intermediate_GenFct} so that $(\mu_{kl}-\mathfrak{A})^{j-1}\phi_{klj}$ will be \eqref{EigenFct_Order1}, or equivalently, $(\mu_{kl}-\mathfrak{A})^j\phi_{klj} = 0$. It remains to show that $\phi_{klj}\in D(\mathfrak{A})$. This condition is equivalent to 
\begin{align*}
    \lambda_0(0)\phi_{klj}(0) - A_d\lambda_0(1)\phi_{klj}(1) = 0.
\end{align*}
By considering \eqref{Gen_EigFct} and $\lambda_0(0)\neq 0$, the previous equation can be rewritten as
\begin{align*}
\omega_j - A_de^{-\mu_{kl}\eta(1)}\left[\omega_j + \sum_{m=0}^{j-2}\omega_{j-1-m}\Omega_{m+1}(1)\right] = 0,
\end{align*}
which is equivalent to
\begin{align}
(\rho_k-A_d)\omega_j - A_d\sum_{m=0}^{j-2}\omega_{j-1-m}\Omega_{m+1}(1) = 0.\label{EquaD(A)}
\end{align}
According to the construction of the vectors $\omega_m$, see \eqref{Construction_omega}, the equation \eqref{EquaD(A)} is satisfied, that is, $\phi_{klj}\in D(\mathfrak{A})$. Assume now that $A_d$ has a basis of eigenvectors denoted by $\{v_k\}_{1\leq k\leq n}$. To show that the set $\{\phi_{kl}\}_{1\leq k\leq n, l\in\mathbb{Z}}$ is complete, we take $f\in X$ and we assume that $\langle \phi_{kl},f\rangle_X = 0$ for all $1\leq k\leq n$ and for all $l\in\mathbb{Z}$. The equation $\langle\phi_{kl},f\rangle_X = 0$ is equivalently written as
\begin{equation}
\lambda_0(0)\int_0^1 f^*(\zeta)v_k \lambda_0(\zeta)^{-1}e^{-\mu_{kl}\int_0^\zeta\lambda_0(\tau)^{-1}\d\tau}\d\zeta = 0.
\label{Orthogonality}
\end{equation}
By using the change of variables $\tilde{\eta} = q(\zeta) := \eta(\zeta)\eta(1)^{-1}$, \eqref{Orthogonality} is equivalent to 
\begin{align*}
\int_0^1 f^*(q^{-1}(\tilde{\eta}))v_k e^{\tilde{\eta}\left(\log\vert\rho_k\vert + i\theta_k\right)} e^{2\pi i l\tilde{\eta}}\d\tilde{\eta} = 0.
\end{align*}
Note that by the positivity of $\lambda_0$, $q$ is a bijection from $[0,1]$ to $[0,1]$, and so $q^{-1}$ is well-defined. By using the fact that $\{e^{2\pi i l\zeta}\}_{l\in\mathbb{Z}}$ is an orthonormal basis of $\L^2(0,1;\mathbb{C})$, we get that 
\begin{align*}
f^*(q^{-1}(\tilde{\eta}))v_k e^{\tilde{\eta}\left(\log\vert\rho_k\vert + i\theta_k\right)} = 0
\end{align*}
in $\L^2(0,1;\mathbb{C})$. Since $\{v_k\}_{1\leq k\leq n}$ is a basis of $\mathbb{C}^n$ and since the function $q$ is a bijection, we have that $f\equiv 0$ in $X$.
Now we consider the inner product \eqref{Weighted_InnerProduct} in which the matrix $\mathcal{W}$ is defined via $v_k^*\mathcal{W}(\zeta)v_l = w_k(\zeta)\delta_{kl}$ with $w_k(\zeta) := \frac{\displaystyle e^{2\eta(1)^{-1}\eta(\zeta)\log\vert\rho_k\vert}}{\vert \lambda_0(0)\vert^2\eta(1)}$. With this definition, taking $1\leq k, \tilde{k}\leq n, k\neq\tilde{k}$ implies that $\langle \phi_{kn},\phi_{\tilde{k}m}\rangle_\mathcal{W} = 0$ for all $m, n\in\mathbb{Z}$. Now we pick any $1\leq k\leq n$ and any $m, n\in\mathbb{Z}$. There holds
\begin{align*}
&\langle \phi_{kn}, \phi_{km}\rangle_\mathcal{W}\\
&= \int_0^1\lambda_0(\zeta)^{-1}v_k^*\mathcal{W}(\zeta)v_k\vert\lambda_0(0)\vert^2e^{-(\overline{\mu_{km}} + \mu_{kn})\int_0^\zeta\lambda_0(\tau)^{-1}\d\tau}\d\zeta\\
&= \int_0^1\lambda_0(\zeta)^{-1}w_k(\zeta)\vert\lambda_0(0)\vert^2e^{-(\overline{\mu_{km}} + \mu_{kn})\int_0^\zeta\lambda_0(\tau)^{-1}\d\tau}\d\zeta.
\end{align*}
By using \eqref{Eigenvalues_A} together with the expression of $w_k$ implies that
\begin{align*}
\langle \phi_{kn}, \phi_{km}\rangle_\mathcal{W} &= \int_0^1 e^{-2\pi i(n-m)\tilde{\eta}}\d\tilde{\eta} = \delta_{nm},
\end{align*}
where we used $\tilde{\eta} := q(\zeta)$ as change of variables.
\endproof

\begin{remark}
\begin{itemize}
\item One of the interesting features of Theorem \ref{Main_Thm} is that we only need to compute a finite number of eigenvalues to generate all the eigenvalues of the operator $\mathfrak{A}$.
\item The eigenvalues of the operator $\mathfrak{A}$ mirror some of the properties of the eigenvalues of $-K^{-1}LP(1)$. As an example, consider $\rho_m = \rho_n, m\neq n$. This implies that $\mu_{ml} = \mu_{nl}$ for all $l\in\mathbb{Z}$. Moreover, if $\rho$ is an eigenvalue of $-K^{-1}LP(1)$ with algebraic or geometric multiplicity larger than 1, so is the corresponding eigenvalue of $\mathfrak{A}$, see \eqref{Eigenvalues_A}. It is also worth noting that the generalized eigenvectors of $\mathfrak{A}$ are directly linked to the generalized eigenvectors of $A_d$, see Theorem \ref{Main_Thm}.
\end{itemize}
\end{remark}

\begin{cor}\label{cor_Stab}
The growth bound of the semigroup generated by $\mathfrak{A}$ is given by 
\begin{equation}
\omega_0  = p(1)^{-1}\max_{k=1, \dots, n} \log(\vert\rho_k\vert).
\label{GrowthBound_A}
\end{equation}
Moreover, the $C_0$-semigroup generated by the operator $\mathfrak{A}$ is exponentially stable if and only if all the eigenvalues of the matrix $-K^{-1}LP(1)$ lie in the interior of the unit circle.
\end{cor}

\proof
To establish \eqref{GrowthBound_A}, we refer to \cite[Theorem 2.8]{JacobKaiserZwart}. Therein, it is specified that for our class of systems, the spectrum determined growth assumption (SDGA) is satisfied when $\mathfrak{A}$ is a discrete Riesz-spectral operator, that is, $\omega_0 = s(\mathfrak{A})$, where $s(\mathfrak{A}) := \sup_{s\in\sigma(\mathfrak{A})} \mathrm{Re}(s)$ is the spectral bound of $\mathfrak{A}$. According to \eqref{Eigenvalues_A} and to the SDGA, it is easy to obtain \eqref{GrowthBound_A}. Hence, exponential stability is equivalent to the condition $\max_{k=1, \dots, n} \log(\vert\rho_k\vert) < 0$. This condition is satisfied if and only if all the eigenvalues of $-K^{-1}LP(1)$ are in the interior of the unit circle.
\endproof

The stability condition obtained in Corollary \ref{cor_Stab} may be compared to the one given in \cite[Chapter 3]{BastinCoron_Book}. Therein, the authors consider the class \eqref{StateEquation_Diag_Uniform}--\eqref{BC_Diag_Uniform} in which the matrix $M$ is set to $0$ and the function $\lambda_0$ is replaced by a diagonal matrix with constant coefficients. The sufficient condition that they obtain to guarantee exponential stability has been derived from a Lyapunov analysis. Later in the same reference, the authors consider a more general case with spatially dependent coefficients. They give sufficient conditions for the system to be exponentially stable by means of matrix inequalities, see \cite[Proposition 5.1]{BastinCoron_Book}. When the coefficients are constant, one of the two matrix inequalities turns out to be a condition like the one we have in Corollary \ref{cor_Stab}. However, their conditions remain sufficient and not equivalent.

\section{APPLICATION}\label{Section_Application}

\begin{figure}
\begin{center}
\includegraphics[scale=1]{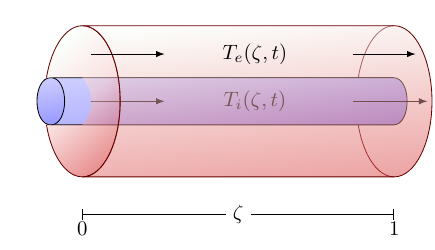}
\end{center}
\caption{Schematic profile view of a co-current heat-exchanger.\label{fig:HE}}
\end{figure}

We consider a co-current heat exchanger of length $1$ as an example to illustrate our main result. A schematic profile view of such a device is depicted in Figure \ref{fig:HE}.
The device consists in two tubes, one immersed in the other. A fluid flows in each tube and the heat exchange occurs only between the two tubes, meaning that we neglect the heat losses with the external environment. The temperature in the internal and the external tubes at position $\zeta\in [0,1]$ and time $t\geq 0$ are denoted by $T_i(\zeta,t)$ and $T_e(\zeta,t)$, respectively. By defining $\mathbf{T}(\zeta,t) = (\begin{smallmatrix}T_i(\zeta,t) & T_e(\zeta,t)\end{smallmatrix})^T$, the PDE governing this system is given by
\begin{equation}
\frac{\partial \mathbf{T}}{\partial t}(\zeta,t) = -v(\zeta)\frac{\partial\mathbf{T}}{\partial \zeta}(\zeta,t) + \left(\begin{smallmatrix}-\alpha_1(\zeta) & \alpha_1(\zeta)\\ \alpha_2(\zeta) & -\alpha_2(\zeta)\end{smallmatrix}\right)\mathbf{T}(\zeta,t),
\label{PDE_Temp}
\end{equation}
where $\alpha_i(\zeta) > 0, i=1, 2, \forall \zeta\in [0,1]$, are heat transfer functions and $v(\zeta) > 0, \forall \zeta\in [0,1]$, is the velocity at which the fluid is flowing, at position $\zeta$. In what follows, we denote by $M(\zeta)$ the matrix $ \left(\begin{smallmatrix}-\alpha_1(\zeta) & \alpha_1(\zeta)\\ \alpha_2(\zeta) & -\alpha_2(\zeta)\end{smallmatrix}\right)$. The initial temperature profile is given by $\mathbf{T}(\zeta,0) = \mathbf{T}_0(\zeta)$. Such a setting may be found for instance in \cite{Maidi_2010}. It is assumed here that the inlet temperature in the internal tube is proportional to the difference temperature between the internal and the external tubes at the outlet, resulting in the following boundary condition $T_i(0,t) = \kappa(T_e(1,t)-T_i(1,t))$, where $\kappa > 0$. In addition, the difference temperature between the outlet and the inlet of the external tube is assumed to be $0$, which is written as $0 = T_e(1,t) - T_e(0,t)$. Those two boundary conditions may then be written in a compact form as
\begin{equation}
0 = \left(\begin{smallmatrix}1 & 0\\ 0 & -1\end{smallmatrix}\right)\mathbf{T}(0,t) + \left(\begin{smallmatrix}\kappa & -\kappa\\ 0 & 1\end{smallmatrix}\right)\mathbf{T}(1,t).
\label{BC_Input}
\end{equation}
The question for this application is to characterize exponential stability of \eqref{PDE_Temp} with \eqref{BC_Input}. In order to put \eqref{PDE_Temp} with \eqref{BC_Input} in the form \eqref{StateEquation_Diag_Uniform}--\eqref{BC_Diag_Uniform}, we introduce the variable $\mathcal{T}(\zeta,t) := v(\zeta)^{-1}\mathbf{T}(\zeta,t)$. This allows to write \eqref{PDE_Temp}, \eqref{BC_Input} equivalently as 
\begin{align}
&\frac{\partial \mathcal{T}}{\partial t}(\zeta,t) = -\frac{\partial}{\partial \zeta}(v(\zeta)\mathcal{T}(\zeta,t)) + M(\zeta)\mathcal{T}(\zeta,t),\label{State_Temp_Change}\\
&0 = -v(0)\left(\begin{smallmatrix}-1 & 0\\ 0 & 1\end{smallmatrix}\right)\mathcal{T}(0,t) - v(1)\left(\begin{smallmatrix}-\kappa & \kappa\\ 0 & -1\end{smallmatrix}\right)\mathcal{T}(1,t).\label{BC_Temp_Change}
\end{align}
Now we compute the solution of the matrix differential equation \eqref{ODE_Matrix}, that is, the solution of
\begin{equation*}
P'(\zeta) = \frac{1}{v(\zeta)}\left(\begin{smallmatrix}-\alpha_1(\zeta) & \alpha_1(\zeta)\\ \alpha_2(\zeta) & -\alpha_2(\zeta)\end{smallmatrix}\right)P(\zeta),\hspace{1cm} P(0) = I.
\end{equation*}
The solution to the previous equation is given by $P(\zeta) := \left(\begin{smallmatrix}P_1(\zeta) & P_2(\zeta) - e^{h(\zeta)}\\ P_1(\zeta)-e^{h(\zeta)} & P_2(\zeta)\end{smallmatrix}\right)$ where 
\begin{align*}
P_i(\zeta) &= 1 - \int_0^\zeta e^{h(\eta)}\alpha_i(\eta)v(\eta)^{-1}\d\eta,\hspace{1cm} i = 1, 2,
\end{align*}
with
\begin{equation*}
h:[0,1]\to\mathbb{R}^-, h(\zeta) := -\int_0^\zeta (\alpha_1(\tau) + \alpha_2(\tau))v(\tau)^{-1}\d\tau.
\end{equation*}
The matrix $P$ allows to transform \eqref{State_Temp_Change}--\eqref{BC_Temp_Change} into the following PDEs
\begin{align}
\frac{\partial z}{\partial t}(\zeta,t) &= -\frac{\partial}{\partial \zeta}(v(\zeta)z(\zeta,t)),\label{State_Temp_Change_P}\\
&0 = -v(0)\left(\begin{smallmatrix}-1 & 0\\ 0 & 1\end{smallmatrix}\right)z(0,t)\nonumber\\
&\hspace{1cm} - v(1)\left(\begin{smallmatrix}-\kappa & \kappa\\ 0 & -1\end{smallmatrix}\right)P(1)z(1,t),\label{BC_Temp_Change_P}
\end{align}
where the variables $z$ and $\mathcal{T}$ are related as $z(\zeta,t) = P(\zeta)^{-1}\mathcal{T}(\zeta,t)$. Note that \eqref{State_Temp_Change_P}--\eqref{BC_Temp_Change_P} is in the same form as \eqref{StateEquation_Diag_Uniform_Chg_Var}--\eqref{BC_Diag_Uniform_Chg_Var} with the matrices $K$ and $L$ given by $K = \left(\begin{smallmatrix}-1 & 0\\ 0 & 1\end{smallmatrix}\right)$ and $L = \left(\begin{smallmatrix}-\kappa & \kappa\\ 0 & -1\end{smallmatrix}\right)$, respectively. Now we concentrate on the well-posedness of \eqref{State_Temp_Change_P}--\eqref{BC_Temp_Change_P}. According to Lemma \ref{Lemma_Well-Posed}, the operator dynamics in \eqref{State_Temp_Change_P}--\eqref{BC_Temp_Change_P} generates a $C_0$-semigroup if and only if the matrix $K$ is invertible. This is trivially satisfied. For that $C_0$-semigroup to be a $C_0$-group, the matrix $L$ needs to be invertible as well. This will be the case if and only if the parameter $\kappa\neq 0$, which is satisfied by assumption on the system parameters. In order to analyze exponential stability, we need to compute the matrix $-K^{-1}LP(1)$. The latter is given by 
\begin{align*}
-K^{-1}LP(1) &= \left(\begin{smallmatrix}-1 & 0\\ 0 & 1\end{smallmatrix}\right)\left(\begin{smallmatrix}-\kappa & \kappa\\ 0 & -1\end{smallmatrix}\right)\left(\begin{smallmatrix}P_1(1) & P_2(1) - e^{h(1)}\\ P_1(1)-e^{h(1)} & P_2(1)\end{smallmatrix}\right)\\
&= \left(\begin{smallmatrix}-\kappa e^{h(1)} & \kappa e^{h(1)}\\ 1 - e^{h(1)} - \int_0^1 e^{h(\eta)}\frac{\alpha_1(\eta)}{v(\eta)}\d\eta & 1 - \int_0^1 e^{h(\eta)}\frac{\alpha_2(\eta)}{v(\eta)}\d\eta\end{smallmatrix}\right).
\end{align*}
The eigenvalues of this matrix are given by
\begin{align*}
\lambda_i &= 2^{-1}\left(1 - \kappa e^{h(1)} - \displaystyle\int_0^1 e^{h(\eta)}\frac{\alpha_2(\eta)}{v(\eta)}\d\eta + (-1)^{i+1} \sqrt{\rho}\right),
\end{align*}
where $\rho = \left[-1+\kappa e^{h(1)} + \int_0^1 e^{h(\eta)}\frac{\alpha_2(\eta)}{v(\eta)}\d\eta\right]^2 + 4\kappa e^{h(1)}$. From the expression of $\rho$, we deduce easily that $\lambda_1>0$ and that $\lambda_2<0$. Observe now that the condition $\lambda_1<1$ holds if and only if 
\begin{align*}
&\left[-1+\kappa e^{h(1)} + \int_0^1 e^{h(\eta)}\frac{\alpha_2(\eta)}{v(\eta)}\d\eta\right]^2 + 4\kappa e^{h(1)}\\
&< \left(1 + \kappa e^{h(1)} + \int_0^1 e^{h(\eta)\frac{\alpha_2(\eta)}{v(\eta)}}\d\eta\right)^2.
\end{align*}
It can then be shown that this condition is satisfied for all $\kappa > 0$. Similarly, the condition $-1<\lambda_2$ is equivalent to 
\begin{align*}
\sqrt{\rho}< 3 - \kappa e^{h(1)} - \displaystyle\int_0^1 e^{h(\eta)}\alpha_2(\eta)v(\eta)^{-1}\d\eta.
\end{align*}
This inequality is valid if and only if 
\begin{align*}
2\kappa e^{h(1)} + \int_0^1 e^{h(\eta)}\alpha_2(\eta)v(\eta)^{-1}\d\eta < 2,
\end{align*}
or equivalently
\begin{equation}
\kappa < 2^{-1}e^{-h(1)}\left(2 - \displaystyle\int_0^1 e^{h(\eta)}\alpha_2(\eta)v(\eta)^{-1}d\eta\right).
\label{ConditionKappa}
\end{equation}
This means that the operator dynamics associated to \eqref{State_Temp_Change_P}--\eqref{BC_Temp_Change_P} generates an exponentially stable $C_0$-semigroup if and only if the parameter $\kappa$ satisfies \eqref{ConditionKappa}.

\section{CONCLUSION AND PERSPECTIVES}
Spectral properties of PDEs representing networks of waves have been studied in this work. In particular, necessary and sufficient conditions guaranteeing well-posedness have been given in Section \ref{Section_Well-Posed}. In Section \ref{Section_Spectral-Analysis}, it has been shown that the operator dynamics being the generator of a $C_0$-group is equivalent for that operator to be a discrete Riesz-spectral operator. In such a case, explicit analytical expressions of the eigenvalues and the eigenvectors have been given, leading to easily checkable stability conditions. Our main results have been applied to characterize exponential stability of a co-current heat exchanger in Section \ref{Section_Application}.

The proposed study in our note opens the door for control design for systems like \eqref{StateEquation_Diag_Uniform}--\eqref{BC_Diag_Uniform}. A first control approach that could be tackled is the positive stabilization as described in \cite{AboAch22}. Therein, a key assumption is that the operator dynamics is a Riesz-spectral operator. Based on this, the authors design a control law that stabilizes the system while keeping the state trajectory positive. This would be an important and a meaningful control problem for \eqref{StateEquation_Diag_Uniform}--\eqref{BC_Diag_Uniform} since the state components for the applications described in the introduction are often required to be positive. As another control method, we mention the Linear-Quadratic optimal control problem tackled in \cite{HastirJacobZwart_LQ}. It is worth noting that the closed-loop system in this reference is again of the form \eqref{StateEquation_Diag_Uniform}--\eqref{BC_Diag_Uniform}. Hence, its operator dynamics could be a discrete Riesz-spectral operator under some matrix condition, which eases a lot the analysis of the closed-loop system in that case. More generally, we refer to \cite{CurtainZwart2020} in which a lot of properties are easier to analyze when the operator under condiseration is a Riesz-spectral operator. By properties, we mean stability, controllability, observability, design of stabilizing compensators to mention a few.

Future works aim at extending our main results to a class of systems similar to \eqref{StateEquation_Diag_Uniform}--\eqref{BC_Diag_Uniform} in which the function $\lambda_0$ could be replaced by a diagonal matrix with possibly different entries. An explicit computation of eigenvalues and eigenfunctions could be worth to investigate in such cases.

\bibliographystyle{plain}
\bibliography{biblio_Main}

\begin{thebibliography}{10}

\bibitem{AboAch22}
B.~Abouza{\"i}d, M.~E. Achhab, J.~N. Dehaye, A.~Hastir, and J.~J. Winkin.
\newblock Locally positive stabilization of infinite-dimensional linear systems
  by state feedback.
\newblock {\em European Journal of Control}, 63:1--13, 2022.

\bibitem{BastinCoron_Book}
G.~Bastin and J.M. Coron.
\newblock {\em Stability and Boundary Stabilization of 1-D Hyperbolic Systems}.
\newblock Progress in Nonlinear Differential Equations and Their Applications.
  Springer International Publishing, 2016.

\bibitem{Curtain_Riesz}
R.~F. Curtain.
\newblock Spectral systems.
\newblock {\em International Journal of Control}, 39(4):657--666, 1984.

\bibitem{CurtainZwart1995}
R.F. Curtain and H.J. Zwart.
\newblock {\em An Introduction to Infinite-Dimensional Linear Systems Theory}.
\newblock Springer-Verlag, springer edition, 1995.

\bibitem{CurtainZwart2020}
R.F. Curtain and H.J. Zwart.
\newblock {\em Introduction to Infinite-Dimensional Systems Theory: A
  State-Space Approach}, volume~71 of {\em Texts in Applied Mathematics book
  series}.
\newblock Springer New York, United States, 2020.

\bibitem{Dunford_PartIII}
N.~Dunford and J.T. Schwartz.
\newblock {\em Linear Operators, Part III: Spectral Operators}.
\newblock Linear Operators. John Wiley \& Sons, New York, 1971.

\bibitem{Guo_SIAM_2001}
B.-Z. Guo.
\newblock Riesz basis approach to the stabilization of a flexible beam with a
  tip mass.
\newblock {\em SIAM J. Control Optim.}, 39(6):1736--1747, 2001.

\bibitem{Guo_SIAM_2002}
B.-Z. Guo.
\newblock Riesz basis property and exponential stability of controlled
  {E}uler--{B}ernoulli beam equations with variable coefficients.
\newblock {\em SIAM J. Control Optim.}, 40(6):1905--1923, 2002.

\bibitem{Guo2019_Book}
B.-Z. Guo and J.-M. Wang.
\newblock {\em Control of Wave and Beam PDEs: The Riesz Basis Approach}.
\newblock Springer International Publishing, Cham, 2019.

\bibitem{Guo_2004_SCL}
B.-Z. Guo and G.-Q. Xu.
\newblock Riesz bases and exact controllability of {$C_0$}-groups with
  one-dimensional input operators.
\newblock {\em Systems \& Control Letters}, 52(3):221--232, 2004.

\bibitem{GuoZwart_2001Riesz}
B.-Z. Guo and H.~J. Zwart.
\newblock {\em Riesz spectral systems}.
\newblock Number 1591 in Memorandum. Department of Applied Mathematics,
  University of Twente, 2001.

\bibitem{HastirJacobZwart_LQ}
A.~Hastir, B.~Jacob, and H.~J. Zwart.
\newblock Linear-{Q}uadratic optimal control for boundary controlled networks
  of waves.
\newblock {\em arXiv:2402.13706}, 2024.

\bibitem{JacobKaiserZwart}
B.~Jacob, J.~T. Kaiser, and H.~J. Zwart.
\newblock Riesz bases of port-{H}amiltonian systems.
\newblock {\em SIAM J. Control Optim.}, 59(6):4646--4665, 2021.

\bibitem{JacobMorrisZwart}
B.~Jacob, K.~Morris, and H.J. Zwart.
\newblock ${C}_0$-semigroups for hyperbolic partial differential equations on a
  one-dimensional spatial domain.
\newblock {\em Journal of Evolution Equations}, 15:493--502, 2015.

\bibitem{JacobZwart}
B.~Jacob and H.J. Zwart.
\newblock {\em Linear Port-{H}amiltonian Systems on Infinite-dimensional
  Spaces}.
\newblock Number 223 in Operator Theory: Advances and Applications. Springer,
  2012.

\bibitem{Luo_Guo}
Z.H. Luo, B.Z. Guo, and {\"O}.~Morg{\"u}l.
\newblock {\em Stability and Stabilization of Infinite Dimensional Systems with
  Applications}.
\newblock Communications and Control Engineering. Springer London, 1999.

\bibitem{Maidi_2010}
A.~Maidi, M.~Diaf, and J.-P. Corriou.
\newblock Boundary control of a parallel-flow heat exchanger by input–output
  linearization.
\newblock {\em Journal of Process Control}, 20(10):1161--1174, 2010.

\bibitem{Mas52}
J.~Massau.
\newblock {\em M\'{e}moire sur l'int\'{e}gration graphique des \'{e}quations
  aux d\'{e}riv\'{e}es partielles}.
\newblock Annales de l’Association des Ing\'{e}nieurs sortis des Ecoles
  sp\'{e}ciales de Gand, 1900-1904.

\bibitem{Hyperbolic_Phillips_1957}
R.~S. Phillips.
\newblock Dissipative hyperbolic systems.
\newblock {\em Transactions of the American Mathematical Society},
  86(1):109--173, 1957.

\bibitem{Hyperbolic_Phillips_1959}
R.~S. Phillips.
\newblock Dissipative operators and hyperbolic systems of partial differential
  equations.
\newblock {\em Transactions of the American Mathematical Society}, 90:193--254,
  1959.

\bibitem{trostorff2023characterisation}
S.~Trostorff and M.~Waurick.
\newblock Characterisation for exponential stability of port-{H}amiltonian
  systems.
\newblock {\em arXiv 2201.10367}, 2023.

\bibitem{Villegas}
J.A. Villegas.
\newblock {\em A Port-{H}amiltonian Approach to Distributed Parameter Systems}.
\newblock PhD thesis, University of Twente, The Netherlands, May 2007.

\bibitem{Guo_Riesz_2003}
G.-Q. Xu and B.-Z. Guo.
\newblock Riesz basis property of evolution equations in {H}ilbert spaces and
  application to a coupled string equation.
\newblock {\em SIAM J. Control Optim.}, 42(3):966--984, 2003.

\bibitem{Zwart_2010}
H.~J. Zwart.
\newblock Riesz basis for strongly continuous groups.
\newblock {\em Journal of Differential Equations}, 249(10):2397--2408, 2010.

\end{thebibliography}

\end{document}